\documentclass[submission,copyright,creativecommons]{eptcs}
\providecommand{\event}{SOS 2007} 

\usepackage{iftex}

\ifpdf
  \usepackage{underscore}         
  \usepackage[T1]{fontenc}        
\else
  \usepackage{breakurl}           
\fi
\usepackage{graphicx}
\usepackage{amsmath}
\usepackage{amsthm}
\title{Using GXWeb for Theorem Proving\\ and Mathematical Modelling}
\author{Philip Todd
\institute{Saltire Software\\ Portland OR USA}
\email{philt@saltire.com}
\and
Danny Aley 
\institute{Saltire Software\\ Portland OR USA}
\email{\quad dannya@saltire.com }
}
\def\titlerunning{Using GXWeb}
\def\authorrunning{P.H. Todd \& D. Aley}
\begin{document}
\maketitle

\begin{abstract}
GXWeb is the free browser based version of the symbolic geometry software Geometry Expressions. We demonstrate its use in an educational setting with examples from theorem proving, mathematical modelling and loci and envelopes.

\end{abstract}
\section{Introduction}
One approach to introducing automated geometrical deduction into an educational environment is to take an existing well accepted Dynamic Geometry System (DGS) and create add-on modules which embody algorithms in automated geometry theorem proving \cite{GeogebraDiscovery}.  An advantage of this approach is that new UI can be added incrementally to a familiar DGS, reducing the overhead for the user of adopting a new tool.  An advantage of using GeoGebra specifically as a platform is that its open source nature allows researchers to focus on the geometry theorem automation, leaving the broader geometry interface to others.

In this presentation, we demonstrate software with a fundamentally different architecture, which achieves many of the same goals, but with its own distinct advantages and disadvantages.  We focus on the use of GXWeb in education in the context of geometry theorem proving, and two problems in applied optics.

\section{GXWeb}
\begin{figure}[h]%
\centering
\includegraphics[width=0.95\textwidth]{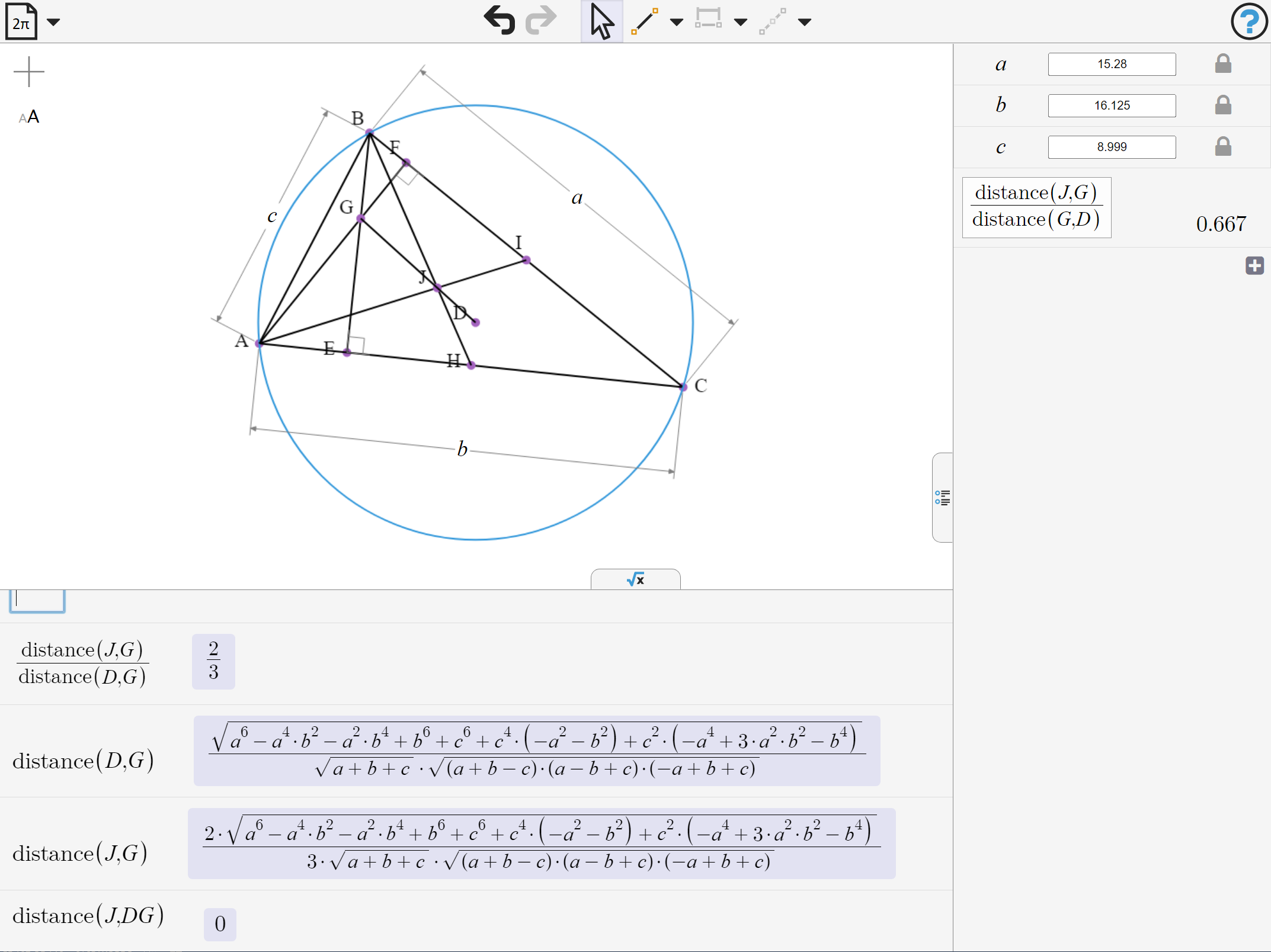}
\caption{The circumcenter $D$, incenter $J$ and orthocenter $G$ of triangle $ABC$ are constructed on the diagram. The numeric panel, to the right shows values for lengths $a$, $b$, $c$, and for the ratio $|JG|/|GD|$ The symbolic panel, below shows exact values for the distance between $J$ and $GD$, along with $|DG|$, $|JG|$ and their ratio.}
\label{fig1}
\end{figure}
GXWeb \cite{GXWeb} is the free browser based version of the symbolic geometry system Geometry Expressions- \cite{todd2020symbolic}.  It maintains both a numeric model of the geometry (similar to that maintained by a typical DGS), but in parallel it also maintains a symbolic model.  User interface is provided (Figure \ref{fig1}) which allows both models to be accessed.  Symbolic inputs are facilitated by a constraint based layer which sits on top of the underlying DGS.  This layer allows distances and angles to be assigned symbolic values in a very natural way.  In parallel to this, a numeric value is maintained for each indeterminate in the symbolic model.  These values control the relation between the symbolic model and the diagram.  Numeric values assigned to variables may be modified in the numeric panel (Figure \ref{fig1}) and output measurements from the numeric model displayed.

The symbolic panel provides access to the symbolic model.  Measurements made in this panel are symbolic rather than numeric. In a theorem proving context, we may use the fact that the distance between a point and a line is symbolically 0 as proof that the point lies on the line.  The symbolic panel can do more than verify relations, however, when non trivial formulas are generated, the user can migrate from exploration into interactive theorem discovery.

To create a locus in GXWeb, one needs to select a point and specify which parameter should vary.  For an envelope one needs instead to select a line or line segment.

In this presentation we demonstrate these features in four educational explorations.  In the first we look at the classical geometry result involving the relative positions of three triangle centers.  We see first how GXWeb can be used in the role of proving a postulated theorem.  We then examine how the symbolic features can be used to suggest further avenues of exploration.
In the second example, we use the numeric model to make hypotheses about tritangent radii for Pythagorean triangles.  These hypotheses can be confirmed in the symbolic panel.
In the third example we look at a geometrical model of a box solar cooker and use GXWeb to help us find the best angle to open the lid.  Finally we pursue an exploration of the caustic curves caused by reflection in a cylinder (the famous coffee cup caustic). 
\section{Examples}
\begin{figure}[h]%
\centering
\includegraphics[width=0.75\textwidth]{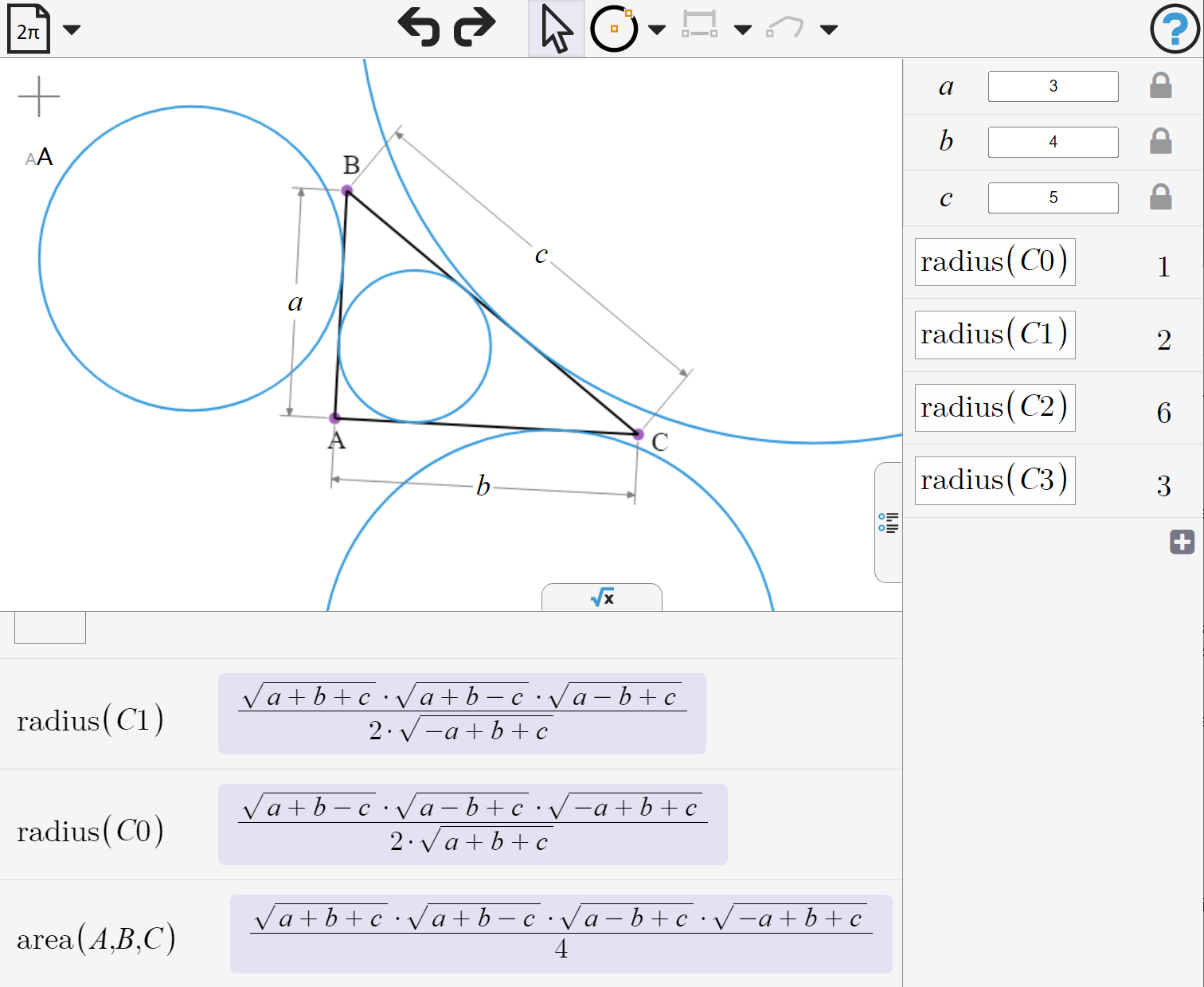}
\caption{Triangle $ABC$ has side lengths $a, b, c$.  The incircle and the three excircles are constructed. Numerical radii are shown for these circles where $a=3$, $b=4$, $c=5$. Symbolic values for two of the radii and for the area of the triangle are shown.}
\label{excircles}
\end{figure}

\subsection{Euler Line}
In figure \ref{fig1}, the sides of a triangle $ABC$ are specified by symbolic lengths $a$, $b$ and $c$, displayed as dimension symbols on the diagram. The numeric values of these variables are displayed in the numeric panel to the right.  The orthocenter is constructed as the intersection between two altitudes, the incenter as the intersection between two medians, and the circumcenter as the center of the constructed circumcircle.  To prove that the three centers are collinear, the user should join two of the points, and inquire, symbolically, for the distance between the third center and this line.  A result of 0 proves the collinearity.

The distances between pairs of centers are complicated formulas involving $a$, $b$ and $c$.  However, their ratio is simple.  Viewed in the symbolic panel, it is a constant.  Viewed in the numeric panel, its value does not vary as the numeric values of the indeterminates vary.

\subsection{Tritangent Circles and Pythagorean Triples}
Numerical experimentation with radii of tritangent circles (Figure \ref{excircles}) leads to a conjecture that for Pythagorean triples, these circles have integer radii.  
Examination of the formulas for the radii gives a starting point for a proof of this conjecture.
It also suggests that multiplying the four radii will provide a significant simplification. This leads to a further conjecture which may be tested numerically, and proved symbolically.

\subsection{Box Solar Cooker}
\begin{figure}[h]%
\centering
\includegraphics[width=0.65\textwidth]{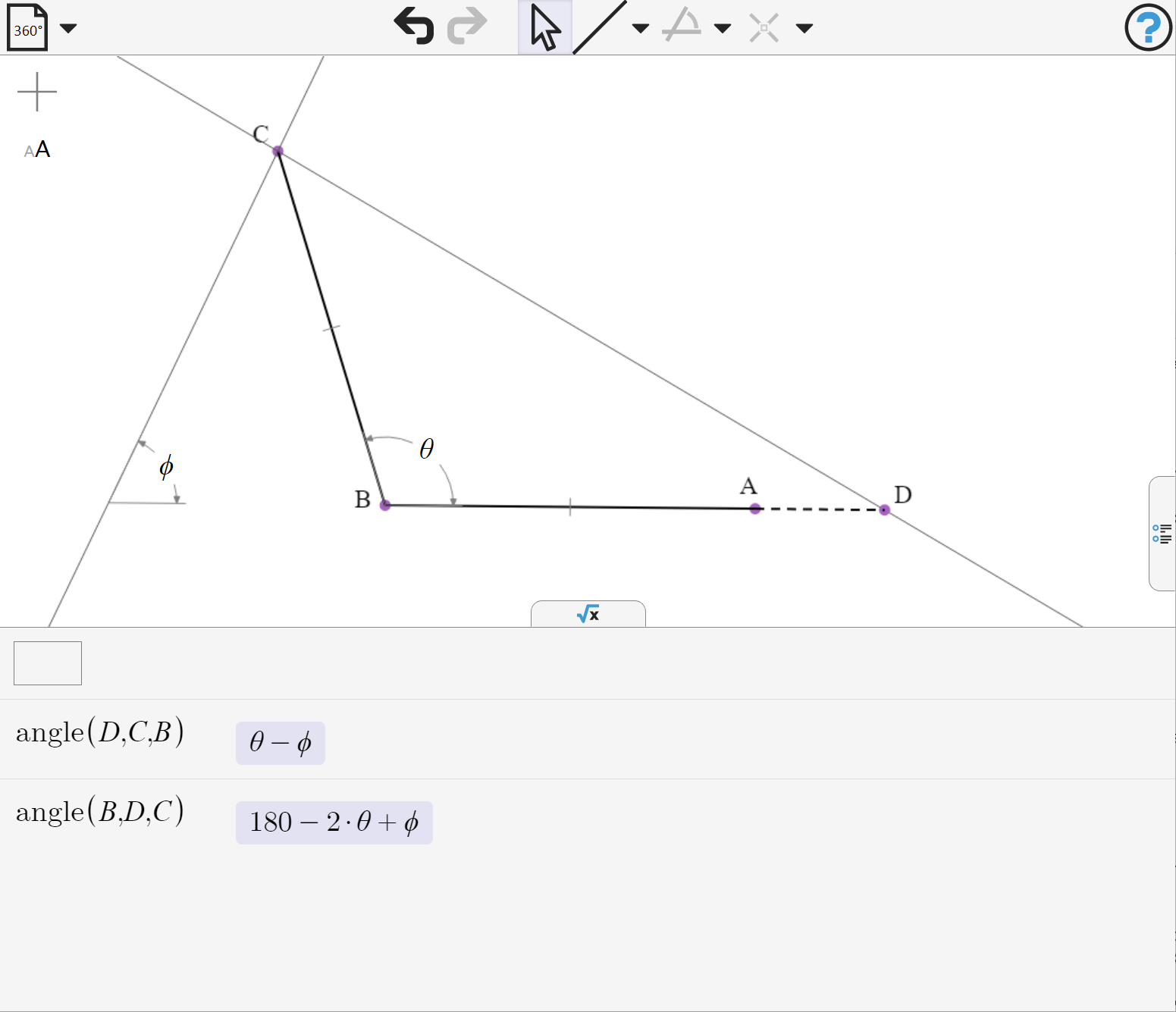}
\caption{A geometric model of a box solar cooker with incident ray at angle $\phi$ to the horizontal and box lid open at angle $\theta$.}
\label{fig2}
\end{figure}
A simple solar cooker is a box with a reflective lid \cite{solarcookers}.  The lid is held open at some angle, and reflected rays are captured in the box.  An interesting theoretical question with this simple apparatus is this:  for given angle of the sun, and assuming the base of the box is on level ground, what is the best angle for the lid.

A geometrical model assumes the sun is at angle $\phi$ to the horizontal and the lid is opened at angle $\theta$. The angle between the ray reflected from the end of the lid to the base can be derived from the model.  Some consideration leads to the postulate that the best angle should reflect light from the outer edge of the lid to the outer edge of the box, thus capturing as many rays as possible, but not spilling any. The angle value generated by GXWeb can be used to set up a linear equation whose solution yields the ideal angle.
\subsection{Coffee Cup Caustic}
\begin{figure}[t]%
\centering
\includegraphics[width=0.6\textwidth]{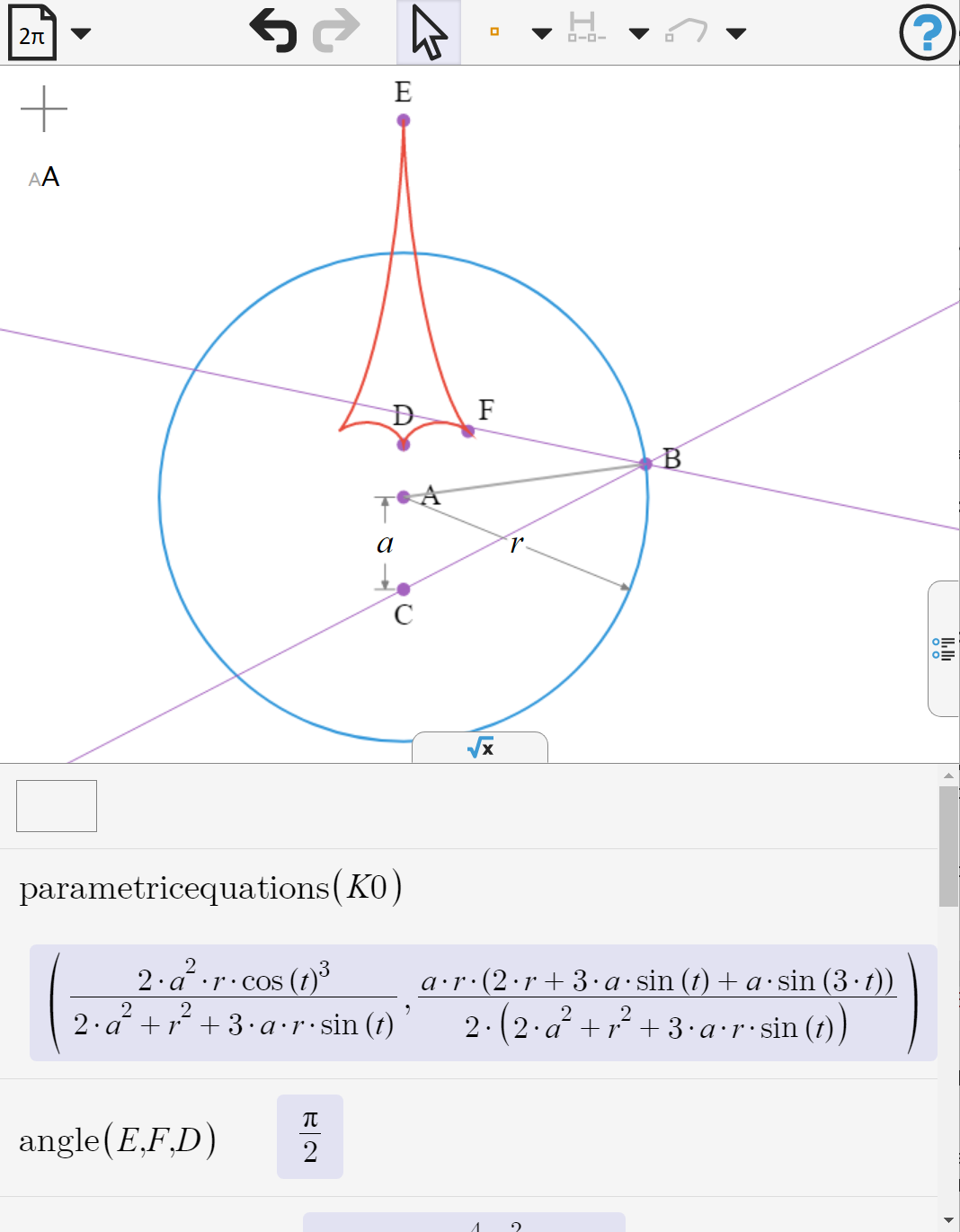}
\caption{The caustic curve due to reflection in a circle of light emanating from point $C$ inside the circle.}
\label{fig3}
\end{figure}
A catacaustic of a curve is the envelope of the reflected rays from some point source (or from a family of parallel rays).  Figure \ref{fig3} shows a geometric construction of this curve. Where the light source lies on the circumference, we can derive an implicit equation for the curve.  For general position of the curve, the parametric equations are more helpful.  By grasping the correspondence between the location on the curve and the parameter of the point of reflection on the circle, students can conjecture the location of points on the cusps of the curve, and discover geometric relations between these points.
\section{Conclusion}

GXWeb has some distinct features which can be advantageous in an educational setting.  The constraint based user interface provides a clean way of setting up and specifying exact geometry problems.  The parallel symbolic and numeric model provide flexibility in moving between a numeric and an algebraic view of a problem.  Loci and envelopes are cleanly defined in terms of a single varying parameter.  Both implicit and parametric equations of the curve may be computed.  Both have their use in different circumstances. 

\nocite{*}
\bibliographystyle{eptcs}
\bibliography{pruned}
\end{document}